\documentclass[12pt]{article}
\usepackage{amsfonts}
\usepackage{latexsym, amssymb, amsmath, amscd, amsfonts, epsfig, graphicx, colordvi, verbatim}
\usepackage{ifpdf}
\usepackage{graphicx}
\usepackage{pstricks}

\newtheorem{thm}{Theorem}[section]

\newtheorem{cor}[thm]{Corollary}

\def\qed{\nopagebreak\hfill{\rule{4pt}{7pt}}
\medbreak}

\def\qed{\nopagebreak\hfill{\rule{4pt}{7pt}}
\medbreak}

\begin{document}

\begin{center}
{\bf \large On  Stanley's Partition Function}

\vskip 3mm
\end{center}

\begin{center}
William Y. C. Chen$^1$, Kathy Q. Ji$^2$, and Albert J. W. Zhu$^3$

 Center for Combinatorics, LPMC-TJKLC\\
   Nankai University\\
    Tianjin 300071, P.R. China

\vskip 1mm

   Email: $^1$chen@nankai.edu.cn, $^2$ji@nankai.edu.cn, $^3$zjw@cfc.nankai.edu.cn

\end{center}

\vskip 6mm \noindent {{\bf Abstract.} Stanley defined a partition
function $t(n)$ as the number of partitions $\lambda$ of $n$ such
that the number of odd parts of $\lambda$ is congruent to the
number of odd parts of the conjugate partition $\lambda'$  modulo $4$. We show that $t(n)$ equals the  number of
partitions of $n$ with an even number of hooks of even length. We derive a closed-form formula for the generating function
for the numbers $p(n)-t(n)$.
 As a consequence, we see that $t(n)$
has the same parity as the ordinary partition function $p(n)$ for
any $n$. A simple combinatorial explanation of this fact is also provided.



\section{Introduction}

This note is concerned with the partition
function $t(n)$ introduced by Stanley \cite{Sta02, Sta05}. We shall give a combinatorial interpretation of $t(n)$ in terms
of hook lengths and shall prove that $t(n)$
and the partition function $p(n)$ have the
same parity. Moreover, we compute the generating function
for $p(n)-t(n)$ and related generating functions.

We shall adopt the common notation on partitions in Andrews
\cite{And76}} or Andrews and Eriksson \cite{AndE04}. A partition
$\lambda =(\lambda_1,\, \lambda_2,\, \lambda_3,\, \cdots,\,
\lambda_r)$ of a nonnegative integer $n$ is a nonincreasing
sequence of nonnegative integers such that the sum of the
components $\lambda_i$ equals $n.$ A part is meant to be a
positive component, and the number of parts of $\lambda$ is called
the length, denoted $l(\lambda)$. The conjugate partition of
$\lambda$ is defined by
$\lambda'=(\lambda'_1,\,\lambda'_2,\,\ldots,\,\lambda'_t),$ where
$\lambda'_i\ (1\le i\le t,\, t=l(\lambda)$) is the number of parts
in $(\lambda_1,\,\lambda_2,\ldots,\,\lambda_r)$ which are greater
than or equal to $i.$ The number of odd parts in
$\lambda=(\lambda_1,\,\lambda_2,\ldots,\,\lambda_r)$ is denoted by
$\mathcal{O}(\lambda).$

For $|q|<1$, the $q$-shifted
factorial is defined by
\begin{equation*}\label{q notation finite}
(a;q)_{n}:=(1-a)(1-aq)\cdots(1-aq^{n-1}),\quad n\ge 1,
\end{equation*}
and
\begin{equation*}\label{q notation infinite}
(a;q)_{\infty}:=(1-a)(1-aq)(1-aq^2)\cdots.
\end{equation*}

Stanley \cite{Sta02, Sta05}   introduced the partition function
$t(n)$ as the number of partitions $\lambda $ of $n$ such that
¡¡$\mathcal{O(\lambda)}\equiv \mathcal{O(\lambda')}\ \mbox{(mod
4)}$, and obtained the following formula
\begin{equation}\label{tpf}
t(n)=\frac{1}{2}\left(p(n)+f(n)\right),
\end{equation}
where $p(n)$ is the number of partitions of $n$ and
\begin{equation}\label{fng}
\sum_{n=0}^{\infty}f(n)q^n=\prod_{i\ge 1}\frac{(1+q^{2i-1})}{(1-q^{4i})(1+q^{4i-2})^2}.
\end{equation}

 Andrews \cite{And04} obtained the following closed-form formula for the
 generating function of  $t(n)$
\begin{equation}\label{an-tn}
\sum_{n=0}^{\infty}t(n)q^n
=\frac{(q^2;q^2)_{\infty}^2(q^{16};q^{16})_{\infty}^5}
{(q;q)_{\infty}(q^4;q^4)_{\infty}^5(q^{32};q^{32})_{\infty}^2}.
\end{equation}
He also derived the congruence relation
\begin{equation}
t(5n+4)\equiv 0 \;({\rm mod } \;  5).
\end{equation}

In this note, we shall consider the complementary partition function
of $t(n)$, namely, the partition function $u(n)=p(n)-t(n)$, which
is the number of partitions $\lambda$ of $n$ such that
$\mathcal{O}(\lambda)\not\equiv\mathcal{O}(\lambda')\;({\rm mod }
\; 4)$. We obtain a closed-form formula for the generating
function of $u(n)$ which implies that Stanley's partition function
$t(n)$ and ordinary partition function $p(n)$ have the same parity
for any $n$. We also present a
simple combinatorial explanation of this fact.
Then we derive formulas for the generating functions for the
numbers $u(4n), u(4n+1), u(4n+2)$ and $u(4n+3)$ which are
analogous to the
 formulas  for the
partition function  $t(n)$ due to Andrews \cite{And04}. In the
last section, we find combinatorial interpretations for $t(n)$ and
$u(n)$ in terms of  hooks of even length.

\section{The generating function formula}

 We shall derive a  formula for the partition
function $u(n)=p(n)-t(n)$. The proof is similar to
 Andrews' proof of (\ref{an-tn}) for $t(n)$.
As a consequence,
one sees that $t(n)$ and $p(n)$ have the same parity for any $n$. This fact also has a simple
combinatorial interpretation.
We shall also compute the generating functions for the numbers
$u(4n), u(4n+1), u(4n+2)$ and $u(4n+3)$.

\begin{thm}\label{theorem on u(n)}
 We have
\begin{equation}\label{generating function of u(n)}
\sum_{n=0}^{\infty}u(n)q^n=\frac{2q^2(q^2;q^2)^2_{\infty}(q^8;q^8)^2_{\infty}(q^{32};q^{32})^2_{\infty}}{(q;q)_{\infty}(q^4;q^4)^5_{\infty}(q^{16};q^{16})_{\infty}}
\end{equation}
\end{thm}

{\noindent \it Proof.}
By (\ref{pf2}) and (\ref{fng}), we find
\begin{align*}
\sum_{n=0}^{\infty}u(n)q^n &=\frac{1}{2}\left(\frac{1}{(q;q)_{\infty}}-
\frac{(-q;q^2)_{\infty}}{(q^4;q^4)_{\infty}
(-q^2;q^4)_{\infty}^2}\right)\\[5pt]
&=\frac{1}{2}\left(\frac{(-q;q^2)_{\infty}}{(q^4;q^4)_{\infty}
(q^2;q^4)_{\infty}^2}-\frac{(-q;q^2)_{\infty}}{(q^4;q^4)_{\infty}
(-q^2;q^4)_{\infty}^2}\right)\\[5pt]
&=\frac{(-q;q^2)_{\infty}}{2(q^4;q^4)^2_{\infty}
(q^2;q^4)_{\infty}^2(-q^2;q^4)_{\infty}^2}
\left(  (q^4;q^4)_{\infty}(-q^2;q^4)_{\infty}^2 -
(q^4;q^4)_{\infty}(q^2;q^4)_{\infty}^2 \right).
\end{align*}
Using Jacobi's triple product identity \cite[p.10]{Ber04}
\begin{equation}\label{jtp}
\sum_{n=-\infty}^{\infty} z^n q^{n^2}
=(-zq;q^2)_{\infty} (-q/z;q^2)_{\infty} (q^2;q^2)_{\infty},
\end{equation}
we see that
\begin{equation} \label{q1}
(q^4;q^4)_{\infty} (-q^2;q^4)_{\infty}^2=\sum_{n=-\infty}^{\infty} q^{2n^2}
\end{equation}
and
\begin{equation}\label{q2}
(q^4;q^4)_{\infty}(q^2;q^4)_{\infty}^2=\sum_{n=-\infty}^{\infty}(-1)^n q^{2n^2}.
\end{equation}
Clearly,
\begin{equation}
\sum_{n=-\infty}^{\infty} q^{2n^2}-\sum_{n=-\infty}^{\infty}(-1)^n q^{2n^2}=2\sum_{n=-\infty}^{\infty}q^{2(2n+1)^2}
\end{equation}
It follows that
\begin{align}
\sum_{n=0}^{\infty}u(q)q^n&=\frac{(-q;q^2)_{\infty}}{(q^4;q^4)^2_{\infty}
(q^2;q^4)_{\infty}^2(-q^2;q^4)_{\infty}^2}
\sum_{n=-\infty}^{\infty}q^{2(2n+1)^2}\nonumber \\[5pt]
&=\frac{q^2(-q;q^2)_{\infty}}{(q^4;q^4)^2_{\infty}(q^4;q^8)_{\infty}^2}\sum_{n=-\infty}^{\infty}q^{8n^2+8n}.
\label{un-3}
\end{align}
 Jacobi's triple product identity yields
\begin{equation}
\sum_{n=-\infty}^{\infty}q^{8n^2+8n}=(-q^{16};q^{16})_{\infty}(-1;q^{16})_{\infty}(q^{16};q^{16})_{\infty}.
\end{equation}
Observe that
\begin{equation}
(-1;q^{16})_{\infty}=2(-q^{16};q^{16})_{\infty}.
\end{equation}
In view of (\ref{un-3}), we get
\begin{align*}
\sum_{n=0}^{\infty}u(q)q^n&=\frac{2q^2(-q^{16};q^{16})_{\infty}
(-q^{16};q^{16})_{\infty}(-q;q^2)_{\infty}
(q^{16};q^{16})_{\infty}}{(q^4;q^4)^2_{\infty}(q^4;q^8)_{\infty}^2}\\[5pt]
&=\frac{2q^2(q^{32};q^{32})_{\infty}(-q;q^2)_{\infty}(-q^{16};q^{16})_{\infty}}{(q^4;q^4)^2_{\infty}(q^4;q^8)_{\infty}^2}.
\end{align*}
Now,
\begin{equation}
(-q;q^2)_{\infty}=
\frac{(q^2;q^2)^2_{\infty}}{(q;q)_{\infty}(q^4;q^4)_{\infty}},
\end{equation}
\begin{equation}
(q^4;q^8)_{\infty}=\frac{(q^4;q^4)_{\infty}}{(q^8;q^8)_{\infty}}
\end{equation}
and
\begin{equation}
(-q^{16};q^{16})_{\infty}=\frac{(q^{32};q^{32})_{\infty}}{(q^{16};q^{16})_{\infty}}.
\end{equation}
Consequently,
\begin{align*}
\sum_{n=0}^{\infty}u(q)q^n&=\frac{2q^2(q^{32};q^{32})_{\infty}
(q^8;q^8)^2_{\infty}(q^2;q^2)^2_{\infty}
(q^{32};q^{32})_{\infty}}{(q^4;q^4)^2_{\infty}
(q^4;q^4)_{\infty}^2(q;q)_{\infty}(q^4;q^4)_{\infty}
(q^{16};q^{16})_{\infty}}\\[5pt]
&=\frac{2q^2(q^2;q^2)^2_{\infty}
(q^8;q^8)^2_{\infty}(q^{32};q^{32})^2_{\infty}}
{(q;q)_{\infty}(q^4;q^4)^5_{\infty}(q^{16};q^{16})_{\infty}}.
\end{align*}
This completes the proof. \qed

\begin{cor} \label{Equidistriubtion}
For   $n\geq 0$,
\[ t(n)\equiv p(n) \; ({\rm mod } \; 2).\]
\end{cor}

 We remark that there is a  simple combinatorial explanation of the above parity property. First, we observe that for any
partition $\lambda$ of $n$,
\begin{equation}\label{on}
 O(\lambda) \equiv O(\lambda') \pmod
2
 \end{equation}
 because we have both  $O(\lambda) \equiv n \pmod
2$ and $O(\lambda') \equiv n \pmod 2$. By the definition of $u(n)$
and the relation (\ref{on}), we see that $u(n)$ equals the number
of partitions of $n$ such that \begin{equation}\label{on2}
O(\lambda) -O(\lambda') \equiv 2 \pmod 4.
\end{equation}
Suppose $\lambda$ is a partition counted by $u(n)$. From
(\ref{on2}) it is evident that its conjugation $\lambda'$ is also
counted by $u(n).$ Once more, from (\ref{on2}) we deduce that
$O(\lambda)$ and $O(\lambda')$ are not equal, so that $\lambda$ is
different from $\lambda'$. Thus we reach the conclusion that
$u(n)$ must be even, and so $t(n)$ has the same parity as $p(n)$
since $p(n)=t(n)+u(n)$.

 From (\ref{tpf}) it follows that
 \begin{equation}\label{pf2}
 u(n) = p(n)-t(n)={p(n)-f(n)\over 2}.
  \end{equation}
  So we have the following congruence relation.

\begin{cor} For $n\geq 0$,
 \[ f(n)\equiv p(n) \;({\rm mod}\; 4).\]
 \end{cor}

Theorem \ref{theorem on u(n)} enables us
to derive the generating functions for $u(4n+i)$, where $i=0,1,2,3$.
Andrews \cite{And04} has obtained formulas for the generating functions
of $t(4n+i)$ for $i=0,1,2,3$.

\begin{thm} We have
\begin{align*}
\sum_{n=0}^{\infty}u(4n)q^{n}&=2q^2(q^{16};q^{16})_{\infty}
(-q;q^{16})_{\infty}(-q^{15};q^{16})_{\infty}V(q)\\[5pt]
\sum_{n=0}^{\infty}u(4n+1)q^{n}&=2q(q^{16};q^{16})_{\infty}
(-q^3;q^{16})_{\infty}(-q^{13};q^{16})_{\infty}V(q)\\[5pt]
\sum_{n=0}^{\infty}u(4n+2)q^{n}&=2(q^{16};q^{16})_{\infty}
(-q^{7};q^{16})_{\infty}(-q^{9};q^{16})_{\infty}V(q)\\[5pt]
\sum_{n=0}^{\infty}u(4n+3)q^{n}&=2(q^{16};q^{16})_{\infty}
(-q^{5};q^{16})_{\infty}(-q^{11};q^{16})_{\infty}V(q)
\end{align*}
where
$$V(q)=\frac{(q^2;q^2)^2_{\infty}(q^8;q^8)^2_{\infty}}{(q;q)^5_{\infty}(q^4;q^4)_{\infty}}.$$
\end{thm}

{\noindent \it Proof.} By Theorem \ref{theorem on u(n)}, we find
\begin{align*}
\sum_{n=0}^{\infty}u(n)q^n&=
\frac{2q^2(q^2;q^2)^2_{\infty}}{(q;q)_{\infty}}V(q^4)\\[5pt]
&=\frac{2q^2(q^2;q^2)_{\infty}}{(q;q^2)_{\infty}}V(q^4)
\end{align*}
Since
\begin{equation}
\frac{1}{(q;q^2)_{\infty}}=(-q;q)_{\infty}
\end{equation}
and
\begin{equation}
(q^2;q^2)_{\infty}=(q;q)_{\infty}(-q;q)_{\infty},
\end{equation}
we have
\begin{align*}
\sum_{n=0}^{\infty}u(n)q^n&=2q^2(q;q)_{\infty}(-q;q)_{\infty}(-q;q)_{\infty}V(q^4)\\
&=q^2(q;q)_{\infty}(-1;q)_{\infty}(-q;q)_{\infty}V(q^4).
\end{align*}
Using Jacobi's triple product identity, we get
\begin{equation}
(q;q)_{\infty}(-1;q)_{\infty}(-q;q)_{\infty}=\sum_{n=-\infty}^{\infty}q^{\frac{n(n+1)}{2}}.
\end{equation}
Thus we have
\begin{equation}\label{un-5}
\sum_{n=0}^{\infty}u(n)q^n=q^2\sum_{n=-\infty}^{\infty}
q^{\frac{n(n+1)}{2}}V(q^4)
=2q^2\sum_{n=0}^{\infty}q^{\frac{n(n+1)}{2}}V(q^4).
\end{equation}
It is easy to check that
\begin{equation}
\sum_{n=0}^{\infty}q^{\frac{n(n+1)}{2}}
=\sum_{n=-\infty}^{\infty}q^{2n^2-n}.
\end{equation}
In view of (\ref{un-5}), we get
\begin{align}
\sum_{n=0}^{\infty}u(n)q^n&=2q^2
\sum_{n=-\infty}^{\infty}q^{2n^2-n}V(q^4)\nonumber \\[5pt]
&=2q^2\sum_{i=0}^{3}\sum_{k=-\infty}^{\infty}q^{2(4k+i)^2-(4k+i)}V(q^4).
\label{un-7}
\end{align}
For $i=0,$ extracting the terms of the form $q^{4j+2}$ in (\ref{un-7})
for $j=0, 1,2,\cdots$,
we obtain
$$\sum_{n=0}^{\infty}u(4n+2)q^{4n+2}=2q^2\sum_{j=-\infty}^{\infty}q^{32j^2-4j}V(q^4).$$
Again,  Jacobi's triple product identity gives
\begin{equation}
\sum_{j=-\infty}^{\infty}q^{32j^2-4j}=(q^{64};q^{64})_{\infty}(-q^{28};q^{64})_{\infty}(-q^{36};q^{64})_{\infty}.
\end{equation}
Hence we get
$$\sum_{n=0}^{\infty}u(4n+2)q^{4n+2}
=2q^2(q^{64};q^{64})_{\infty}(-q^{28};q^{64})_{\infty}
(-q^{36};q^{64})_{\infty}V(q^4),$$
which simplifies to
$$\sum_{n=0}^{\infty}u(4n+2)q^{n}=2(q^{16};q^{16})_{\infty}(-q^{7};q^{16})_{\infty}(-q^{9};q^{16})_{\infty}V(q).$$

The remaining cases can be verified using similar arguments. This completes the proof. \qed

\section{Combinatorial interpretations
for $t(n)$ and $u(n)$}

In \cite[Proposition 3.1]{Sta02},  Stanley found three partition  statistics
 that have the same parity as
$(\mathcal{O}(\lambda)-\mathcal{O}(\lambda'))/2$, and gave several
combinatorial interpretations for $t(n)$. We shall present
combinatorial interpretations of partition functions $t(n)$ and
$u(n)$ in terms of the number of hooks of even length. For the definition of hook lengths,
 see Stanley \cite[p. 373]{Sta99}.
A hook of even length is called an even hook.  The following
theorem  shows that the number of even hooks has the same parity
as $(\mathcal{O}(\lambda)-\mathcal{O}(\lambda'))/2.$

\begin{thm}\label{3main}
For any partition $\lambda$ of $n$, $\mathcal{O}(\lambda) \equiv
 \mathcal{O}(\lambda') \pmod{4}$ if and only if $\lambda$ has  an even number of
 even hooks.
 \end{thm}

 \noindent{\it Proof.}  We use induction on  $n$. It is clear that Theorem \ref{3main} holds for $n=1$. Suppose that it is true for all
partitions of $n$. We aim to show that the conclusion also holds for
all partitions of $n+1$. Let $\lambda$ be a partition of $n+1$ and
$v=(i,j) $ be any an inner corner of the Young diagram of $\lambda$, that is, the removal of
 the square $v$ gives a Young diagram of a partition of $n$. Let $\lambda^-$
denote the partition obtained by removing the square $v$ from the
Young diagram of $\lambda$. We use $H_e(\lambda)$ to denote the
number of squares with even hooks in the Young diagram of
$\lambda$. We claim that
\begin{equation}\label{temp}
H_e(\lambda)\equiv H_e(\lambda^-) \pmod 2 \qquad
\Longleftrightarrow \qquad \lambda_i \equiv \lambda'_j \pmod 2.
\end{equation}
Let $\mathcal{T}(\lambda,v)$ denote the set of all squares in the
Young diagram of $\lambda$ which are in the same row as $v$ or in
the same column as $v.$ After removing the square $v$ from the
Young diagram of $\lambda$, the hook lengths of the squares in
$\mathcal{T}(\lambda,v)$ have decreased by one. Meanwhile, the
hook lengths of other squares remain the same. Furthermore, if
$\lambda_i$ and $\lambda'_j$ have the same parity, then the number
of squares in $\mathcal{T}(\lambda,v)$ is even. This implies that
the parity of the number of squares in $\mathcal{T}(\lambda,v)$
of even hook lengths  coincides with the parity of the number of
squares in $\mathcal{T}(\lambda,v)$ with odd hook lengths. Similarly, for the case when $\lambda_i$ and
$\lambda'_j$  have different parities,
it can be shown that  the number of squares in $\mathcal{T}(\lambda,v)$
of even hook length  is of   opposite parity to the number of
squares in $\mathcal{T}(\lambda,v)$ of odd hook length. Hence we arrive at
(\ref{temp}).

By the inductive hypothesis, we see that $\mathcal{O}(\lambda^-)
\equiv
 \mathcal{O}((\lambda^-)') \pmod{4}$ if and only if
 $H_e(\lambda^-)$ is even. For any inner corner $v=(i,j)$ of
 $\lambda$,  if $\lambda_i \equiv \lambda'_j \pmod 2,$ then
$\mathcal{O}(\lambda)\equiv \mathcal{O}(\lambda')\pmod 4$ if and
only if $\mathcal{O}(\lambda^-)\equiv
\mathcal{O}((\lambda^-)')\pmod 4$. By (\ref{temp}), we find that
in this case,  $H_e(\lambda)$ and $H_e(\lambda^-)$ have the same
parity. Thus the assertion holds for any partition $\lambda$ of
$n+1$. The case that
 $\lambda_i \not\equiv \lambda'_j \pmod 2$ can be justified in the same manner. This completes the proof.  \qed

 From Theorem \ref{3main}, we obtain a   combinatorial
interpretation for  Stanley's partition function $t(n)$, which can
be recast as a combinatorial
 interpretation for $u(n)$.

 \begin{thm}\label{parity} The partition function
 $t(n)$ is equal to the number of partitions of $n$ with an even
 number of even hooks, and the partition function $u(n)$ is equal to the
 number of partitions of $n$ with an odd number of even hooks.
 \end{thm}

Combining Theorem \ref{theorem on u(n)} and Theorem \ref{parity}
we have the following consequence.

 \begin{cor}For any $n$, the number of  partitions of $n$ with an odd number of even hooks is always even.
 \end{cor}

Since $f(n)= t(n)-u(n)$,  we see that  $f(n)$ can be interpreted
as  the signed  counting  of partitions of $n$ with respect to the
number of even hooks, as formally stated below.

\begin{thm}
The  function
$f(n)$ equals the number of
partitions of $n$ with an even number of even hooks minus the number of partitions of $n$ with
an odd number of even hooks.
\end{thm}

 \vspace{.2cm}

\noindent{\bf Acknowledgments.} This work was supported by the 973
Project, the PCSIRT Project of the Ministry of Education, and the
National Science Foundation of China.

\end{document}